\documentclass{amsart}

\usepackage{amssymb}
\usepackage{amsmath}
\usepackage{amscd}
\usepackage{graphics}

\theoremstyle{plain}
\newtheorem{thm}{Theorem}[section]

\newtheorem{prob}[thm]{Problem}
\theoremstyle{definition}
\newtheorem{defn}[thm]{Definition}
\newtheorem{exmp}[thm]{Example}

\theoremstyle{remark}
\newtheorem{rem}[thm]{Remark}

\begin{document}

\title{Lantern relations and rational blowdowns}
\author[H.~Endo]{Hisaaki Endo}
\address{Department of Mathematics, Graduate School of Science,
Osaka University, Toyonaka, Osaka 560-0043, Japan}
\email{endo@math.sci.osaka-u.ac.jp}
\author[Y.~Z.~Gurtas]{Yusuf Z. Gurtas}
\address{Department of Mathematics and Computer Science, 
St. Louis University, MO, USA}
\email{ygurtas@slu.edu}
\thanks{The first author is partially supported by Grant-in-Aid for Scientific Research 
(No.18540083), Japan Society for the Promotion of Science.}
\date{August 5, 2008; MSC 2000: primary 57R17, secondary 57N13, 20F38}

\maketitle

\begin{abstract} 
We discuss a connection between the lantern relation in mapping class groups
and the rational blowing down process for $4$-manifolds. 
More precisely, if we change a positive relator 
in Dehn twist generators of the mapping class group 
by using a lantern relation, the corresponding Lefschetz fibration changes into 
its rational blowdown along a copy of the configuration $C_2$. 
We exhibit examples of such rational blowdowns of Lefschetz fibrations whose blowup is 
homeomorphic but not diffeomorphic to the original fibration. 
\end{abstract}

\section{Introduction}

Lefschetz fibrations relate the topology of symplectic $4$-manifolds 
to the combinatorics on positive relators in Dehn twist generators of mapping class groups. 
Fuller introduced a substitution technique for constructing positive relators to obtain 
an example of non-holomorphic Lefschetz fibrations of genus three \cite{smith}, 
\cite{ozbagci}. 
Many constructions of Lefschetz fibrations as positive relators can be interpreted as 
generalizations of his construction (cf. \cite{EN}), while it has been less investigated 
what such substitutions mean geometrically. 

In this paper we study a particular substitution, the lantern substitution (or 
the $L^{\pm 1}$-substitution in short), for positive relators of mapping class groups. 
The corresponding surgical operation on Lefschetz fibrations turns out to be 
the rational blowing down process, which was discovered by Fintushel and Stern \cite{FS}, 
along a copy of the configuration $C_2$ 
(i.e. a $-4$-framed unknot in Kirby diagrams). 
Applying a theorem of Usher \cite{usher}, we give examples of such rational blowdowns of 
Lefschetz fibrations whose blowup is homeomorphic but not diffeomorphic to the original 
fibration. 

In Section 2 we review the lantern relation in mapping class groups 
and define the lantern substitution for positive relators. 
We discuss a relation between lantern relations and rational blowdowns in Section 3 
and state the main theorem in Section 4. 
We then exhibit some examples in Section 5 and end by observing other relations in Section 6. 

The authors are grateful to K. Yasui for helpful comments on the rational blowing down 
process and to N. Monden for drawing beautiful Kirby diagrams in Figure 1 and Figure 2. 

\section{Lantern relations and substitutions}

Let $\Sigma_g$ be a closed oriented surface of genus $g\; (\geq 2)$ and $\mathcal{M}_g$ the 
mapping class group of $\Sigma_g$. 
We denote by $\mathcal{ F}$ the free group 
generated by all isotopy classes $\mathcal{ S}$ of 
simple closed curves on $\Sigma_g$. There is a natural 
epimorphism $\varpi:\mathcal{ F}\rightarrow \mathcal{ M}_g$ which sends 
(the isotopy class of) a simple closed curve $a$ on $\Sigma_g$ 
to the right-handed Dehn twist $t_a$ along $a$. 
We set $\mathcal{ R}:={\rm Ker}\; \varpi$ and call each element of $\mathcal{ R}$ 
a {\it relator} in the generators $\mathcal{ S}$ of $\mathcal{ M}_g$. 
A word in the generators $\mathcal{S}$ is called {\it positive} if 
it includes no negative exponents. 
We put ${}_W(c):=t_{a_r}^{\varepsilon_r}\cdots 
t_{a_1}^{\varepsilon_1}(c)\in \mathcal{ S}$ for $c\in \mathcal{ S}$ and 
$W=a_r^{\varepsilon_r}\cdots a_1^{\varepsilon_1}\in \mathcal{ F} \; 
(a_1, \ldots ,a_r \in \mathcal{ S}, \varepsilon_1,\ldots ,\varepsilon_r\in \{ \pm 1\})$, 
and put ${}_WV:={}_W(c_1)\cdots {}_W(c_s) \in \mathcal{ F}$ for 
$V=c_1 \cdots c_s \in \mathcal{ F}\; (c_1, \ldots c_s \in \mathcal{ S})$. 

We begin with a precise definition of the lantern relation \cite{dehn}, \cite{johnson}. 

\begin{defn}\label{lantern}
Let $a$ and $b$ be simple closed curves on $\Sigma_g$ 
with geometric intersection number $2$ 
and algebraic intersection number $0$. 
We orient $a$ and $b$ locally on a neighborhood of each 
intersection point $p\in a\cap b$ such that 
the intersection number $(a\cdot b)_p$ 
at $p$ is $+1$. Resolving all intersection points according to 
the local orientations, we obtain a new simple closed curve $c$. 
A regular neighborhood of $a\cup b$,
which can be chosen to include $c$, is a genus $0$ subsurface 
$\Sigma$ of $\Sigma_g$ with four boundary components. 
We denote simple closed curves parallel to four boundary components 
of $\Sigma$ by $d_1,d_2,d_3$, and $d_4$. The relation 
\[
t_{d_1}t_{d_2}t_{d_3}t_{d_4}=t_at_bt_c
\]
is called the {\it lantern relation}. We put $L:=L(a,b)=
abcd_4^{-1}d_3^{-1}d_2^{-1}d_1^{-1}\in \mathcal{R}$. 
\end{defn}

Let $\varrho\in\mathcal{R}\, (\varrho\ne 1)$ be a positive relator of $\mathcal{M}_g$. 
Let $a,b,c,d_1,d_2,d_3$, and $d_4$ be curves as in Definition \ref{lantern}. 
Suppose that $\varrho$ includes $d_1d_2d_3d_4$ as a subword: 
$\varrho=U\cdot d_1d_2d_3d_4\cdot V\; (U,V\in \mathcal{F})$. 
Since $\varrho$ and $U\cdot L\cdot U^{-1}$ are both relators of $\mathcal{M}_g$,  
the positive word 
\begin{align*}
\varrho' & =U\cdot abc\cdot V\; 
(=U\cdot abcd_4^{-1}d_3^{-1}d_2^{-1}d_1^{-1}\cdot d_1d_2d_3d_4\cdot V \\
& = U\cdot L\cdot U^{-1}\cdot U\cdot d_1d_2d_3d_4\cdot V 
= U\cdot L\cdot U^{-1}\cdot \varrho\, ) 
\end{align*}
is also a relator of $\mathcal{M}_g$. The length of the word 
$\varrho'$ is equal to that of $\varrho$ minus one. 

\begin{defn} We say that $\varrho'$ is obtained by applying 
an $L$-{\it substitution} to $\varrho$. 
Conversely, $\varrho$ is said to be obtained by applying 
an $L^{-1}$-{\it substitution} to $\varrho'$. 
We also call these two kinds of operations {\it lantern substitutions} (cf. \cite{EN}). 
\end{defn}

We next recall a definition of Lefschetz fibrations (cf. \cite{yukiomat}, \cite{GS}). 

\begin{defn}
Let $M$ be a closed oriented smooth $4$-manifold. A smooth map 
$f:M\rightarrow S^2$ is called a {\it Lefschetz fibration} of 
genus $g$ if it satisfies the following conditions: 

(i) $f$ has finitely many critical values $b_1,\ldots ,b_n\in 
S^2$ and $f$ is a smooth fiber bundle 
over $S^2-\{ b_1,\ldots ,b_n\}$ with fiber $\Sigma_g$; 

(ii) for each $i\; (i=1,\ldots ,n)$, there exists a unique 
critical point $p_i$ in the {\it singular fiber} $f^{-1}(b_i)$ 
such that $f$ is locally written as 
$f(z_1,z_2)=z_1^2+z_2^2$ with respect to some local complex 
coordinates around $p_i$ and $b_i$ which are compatible with 
orientations of $M$ and $S^2$; 

(iii) no fiber contains a $-1$-sphere.

\end{defn}

\begin{rem}\label{LF}
A more general definition can be found in Chapter 8 of \cite{GS}. 
We treat also Lefschetz fibrations with boundary in the proof of Theorem \ref{bd}. 
\end{rem}

Suppose that $g\geq 2$. 
According to theorems of Kas and Matsumoto, 
there exists a one-to-one correspondence between 
the isomorphism classes of Lefschetz fibrations and 
the equivalence classes of positive relators 
modulo simultaneous conjugations 
\[
c_1\cdot \cdots \cdot c_n \sim 
{}_W(c_1)\cdot \cdots \cdot {}_W(c_n), 
\]
and elementary transformations 
\begin{gather*}
c_1\cdot \cdots \cdot c_i\cdot c_{i+1}\cdot \cdots \cdot c_n \sim 
c_1\cdot \cdots \cdot c_{i+1}\cdot {}_{c_{i+1}^{-1}}(c_i)\cdot \cdots \cdot c_n, \\
c_1\cdot \cdots \cdot c_i\cdot c_{i+1}\cdot \cdots \cdot c_n \sim 
c_1\cdot \cdots \cdot {}_{c_i}(c_{i+1})\cdot c_i\cdot \cdots \cdot c_n, 
\end{gather*}
where $c_1\cdots c_n\in \mathcal{R}$ is a positive relator in the generator $\mathcal{S}$ 
and $W\in \mathcal{F}$. 
This correspondence is described by using the holonomy 
(or monodromy) homomorphism 
induced by the classifying map of $f$ restricted on $S^2-\{ b_1,\ldots b_n \}$ 
(cf. \cite{GS}, \cite{yukiomat}, and \cite{endo}). 
We denote (the isomorphism class of) 
a Lefschetz fibration associated to a positive relator $\varrho \in \mathcal{R}$ 
by $M_{\varrho}\rightarrow S^2$. 

\medskip

Let $\varrho,\varrho' \in\mathcal{R}$ be positive relators of $\mathcal{M}_g$ 
and $M_{\varrho},M_{\varrho'}$ the corresponding Lefschetz fibrations over $S^2$, respectively. 
Suppose that the relator $\varrho'$ is obtained by applying 
an $L$-substitution to the relator $\varrho$. 
The Euler characteristic and 
the signature of a Lefschetz fibration $M_{\varrho'}\rightarrow S^2$ 
with monodromy $\varrho'$ are computed as follows: 
\[
e(M_{\varrho'})=e(M_{\varrho})-1, \quad 
\sigma(M_{\varrho'})=\sigma(M_{\varrho})+1
\]
(\cite{EN}, Theorem 4.3 and Proposition 3.12). 
We investigate relations between $M_{\varrho}$ and $M_{\varrho'}$ and 
several properties of them in the subsequent sections.

\section{Rational blowdowns via lantern relations}

Let $\varrho,\varrho' \in\mathcal{R}$ be positive relators of $\mathcal{M}_g$ 
and $M_{\varrho},M_{\varrho'}$ the corresponding Lefschetz fibrations over $S^2$, respectively. 

\begin{thm}\label{bd}
If $\varrho'$ is obtained by applying 
an $L$-substitution to $\varrho$, 
then the $4$-manifold $M_{\varrho'}$ is a rational blowdown of $M_{\varrho}$ 
along a configuration $C_2\subset M_{\varrho}$. 
\end{thm}

\noindent
{\it Proof}. We take a subsurface $\Sigma$ of $\Sigma_g$ and curves 
$a,b,c,d_1,d_2,d_3$, and $d_4$ on $\Sigma$ as in Definition \ref{lantern}. 
Let $N,N'$ be Lefschetz fibrations over $D^2$ with fiber $\Sigma$ corresponding to 
the positive words $d_1d_2d_3d_4, abc$, respectively. 

\smallskip

\scalebox{0.8}{\includegraphics{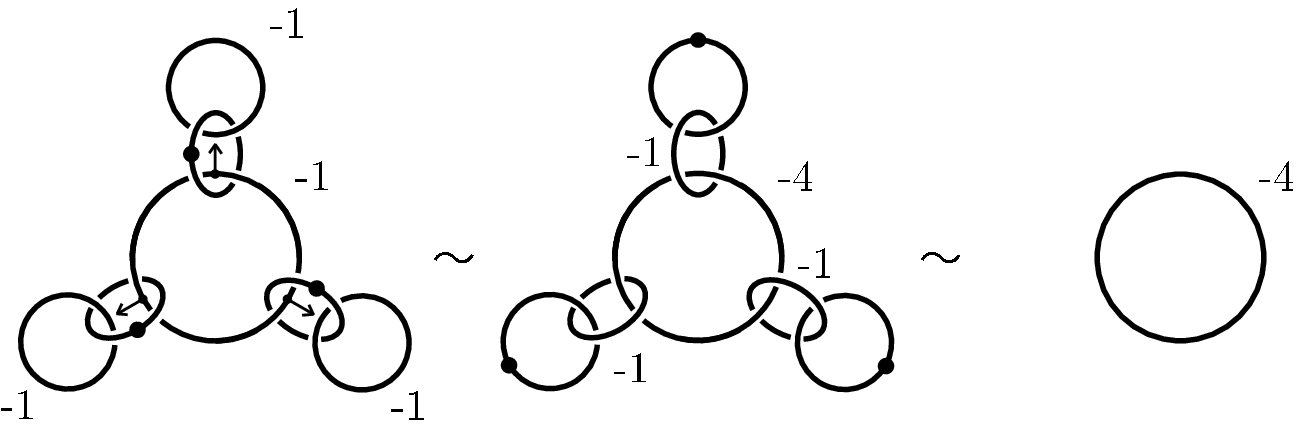}}

\begin{center}
Figure 1
\end{center}

\noindent
Drawing a Kirby diagram of $N$, sliding the central $-1$-framed unknot over 
other three $-1$-framed unknots, and canceling three $1$-handle/$2$-handle pairs, 
we obtain a $-4$-framed unknot (Figure 1). 
Thus $N$ is diffeomorphic to the total space of a $D^2$-bundle over $S^2$ with 
Euler number $-4$, which is denoted by $C_2$ in \cite{FS} (see also \cite{GS}, Section 8.5). 

\scalebox{0.8}{\includegraphics{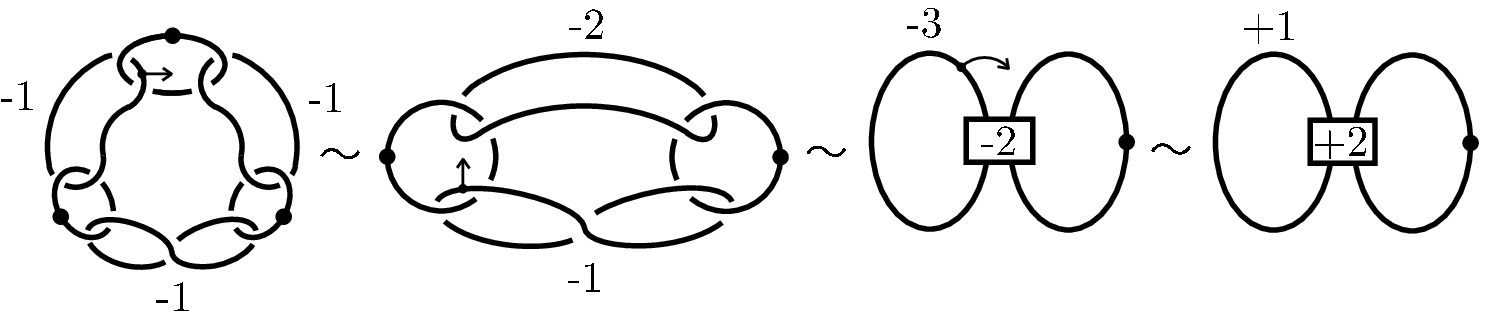}}

\vspace{-3cm}

\begin{center}
Figure 2
\end{center}

\noindent
Drawing a Kirby diagram of $N'$ and sliding and canceling handles as in Figure 2, 
we obtain a pair of a dotted circle and a $+1$-framed unknot with linking number $+2$. 
This means that $N'$ is diffeomorphic to a rational $4$-ball with boundary $L(4,1)$, 
which is denoted by $B_2$ in \cite{FS} (see also \cite{GS}, Section 8.5). 

From construction, 
$N$ (resp. $N'$) can be considered a submanifold of 
$M_{\varrho}$ (resp. $M_{\varrho'}$). 
It is also easily seen that $M_{\varrho}-{\rm int}\, N$ and $M_{\varrho'}-{\rm int}\, N'$ are 
diffeomorphic to each other. Hence we have 
\[
M_{\varrho'}\approx N'\cup_{\partial N}(M_{\varrho}-{\rm int}\, N)
\approx B_2\cup_{L(4,1)}(M_{\varrho}-{\rm int}\, C_2).
\]
This completes the proof of Theorem \ref{bd}. $\square$

\section{Smooth structures}

Let $\varrho,\varrho' \in\mathcal{R}$ be positive relators of $\mathcal{M}_g$ 
and $M_{\varrho},M_{\varrho'}$ the corresponding Lefschetz fibrations over $S^2$, respectively. 
Suppose that $\varrho'$ is obtained by applying $k$ times 
$L$-substitutions ($k\geq 1$), elementary transformations, 
and simultaneous conjugations to $\varrho$. 
Suppose also that $e(M_{\varrho})+\sigma(M_{\varrho})\geq 2$. 
We choose a positive relator $\varsigma\in\mathcal{R}\, (\varsigma\ne 1)$ 
such that $M_{\varsigma}-\nu F$ is simply-connected and either the word $\varsigma$ includes 
at least one separating curve as a factor, or $\sigma(M_{\varrho})+\sigma(M_{\varsigma})$ 
is not divisible by $16$. 
Here $\nu F$ is an open fibered neighborhood of a regular fiber $F$ of $M_{\varsigma}$. 
Taking a fiber sum of $M_{\varrho}$ (resp. $M_{\varrho'}$) with $M_{\varsigma}$, 
we obtain a new Lefschetz fibration 
$M_1:=M_{\varrho}\#_F M_{\varsigma}$ 
(resp. $M_2:=M_{\varrho'}\#_F M_{\varsigma}$) 
with monodromy $\varrho\cdot {}_W\varsigma$ (resp. $\varrho' \cdot {}_{W'}\varsigma$) 
for some $W\in \mathcal{F}$ (resp. $W'\in\mathcal{F}$). 
It is obvious that $\varrho'\cdot {}_{W'}\varsigma$ 
is obtained by applying $k$ times $L$-substitutions, elementary transformations, 
and simultaneous conjugations to $\varrho\cdot {}_W\varsigma$. 

\begin{thm}\label{exotic} The $4$-manifold $M_1$ 
is homeomorphic but not diffeomorphic to a $k$ times blowup 
$M_2\# k\overline{\Bbb{CP}}^2$ of $M_2$. 
Moreover, both of these $4$-manifolds do not dissolve. 
\end{thm}

{\it Proof}. Let $j:F_i\hookrightarrow M_i$ be the inclusion map 
from a general fiber $F_i$ into the total space $M_i\,(i=1,2)$. 
The induced homomorphism 
$j_{\#}:\pi_1(F_i)\rightarrow \pi_1(M_i)$ is surjective and 
the kernel of $j_{\#}$ includes the normal subgroup $N$ of $\pi_1(M_i)$ generated by 
the vanishing cycles of $M_i$ (cf. \cite{ABKP}, Lemma 3.2). 
Since $M_{\varsigma}-\nu F$ is simply-connected and $j_{\#}$ is the composition of 
homomorphisms $\pi_1(F_i)\rightarrow \pi_1(M_{\varsigma}-\nu F)\rightarrow 
\pi_1(M_i)$, the group $\pi_1(M_i)$ must be trivial ($i=1,2$). 

$M_1$ is a non-spin $4$-manifold because either it has a component of 
a separating singular fiber which represents a homology class of square $-1$, 
or $\sigma(M_1)$ is not divisible by $16$. 
It is easily seen from the observation above 
that $e(M_2)=e(M_1)-k$ and $\sigma(M_2)=\sigma(M_1)+k$. 
By virtue of Freedman's classification theorem, both of $M_1$ and 
$M_2\# k\overline{\Bbb{CP}}^2$ is homeomorphic to 
$\# b_2^+(M_1)\Bbb{CP}^2\#b_2^-(M_1)\overline{\Bbb{CP}}^2$ 
because they are simply-connected, non-spin, and 
have the same Euler characteristic and the same signature. 

$M_1$ is a fiber sum $M_{\varrho}\#_F M_{\varsigma}$ of non-trivial Lefschetz fibrations 
$M_{\varrho}$ and $M_{\varsigma}$. 
By Gompf's theorem (\cite{GS}, Theorem 10.2.18), $M_1$ admits a symplectic structure 
with symplectic fibers. 
It follows from a theorem of Usher \cite{usher} that $M_1$ is 
a minimal symplectic $4$-manifold. 
Since $b_2^+(M_1)=b_2^+(M_{\varrho})-b_1(M_{\varrho})+b_2^+(M_{\varsigma})+2g-1>1$, 
$M_1$ does not contain any smooth $-1$-sphere 
as a consequence of Seiberg-Witten theory 
(\cite{taubes1}, \cite{taubes2}, cf. \cite{GS}, Remark 10.2.4(a)). 
On the other hand, $M_2\# k\overline{\Bbb{CP}}^2$ has a natural smooth $-1$-sphere. 
Hence $M_1$ and $M_2\# k\overline{\Bbb{CP}}^2$ can not be diffeomorphic. 

Because $M_1$ and $M_2\# k\overline{\Bbb{CP}}^2$ admit symplectic structure and 
$b_2^+(M_1)>1$, these manifolds can not be 
diffeomorphic to $\# b_2^+(M_1)\Bbb{CP}^2\#b_2^-(M_1)\overline{\Bbb{CP}}^2$ 
(\cite{taubes}, \cite{KMT}, cf. \cite{GS}, Theorem 10.1.14). 
$\square$

\begin{rem} We do not use any explicit property of rational blowdowns 
to prove Theorem \ref{exotic}. The proof above 
is rather similar to that of Theorem 4.8 of \cite{endo}. 
It is likely that $M_{\varrho}$ is homeomorphic but not diffeomorphic to 
$M_{\varrho'}\# k\overline{\Bbb{CP}}^2$ 
(without taking fiber sums with $M_{\varsigma}$) in a general setting. 
On the other hand, a certain rational blowdown along 
$C_2$ happens to be diffeomorphic to an honest blowdown of the original $4$-manifold: 
$E(1)_2\, (\approx E(1)\approx\Bbb{CP}^2\#9\overline{\Bbb{CP}}^2)$ is a rational blowdown of 
$E(1)\# \overline{\Bbb{CP}}^2\,(\approx\Bbb{CP}^2\#10\overline{\Bbb{CP}}^2)$ along $C_2$ 
(\cite{FS}, Proposition 3.2, cf. \cite{GS}, Theorem 8.5.9 and Theorem 8.3.11). 
\end{rem}

\section{Examples}

We apply theorems in previous sections to explicit examples. 

\begin{exmp} 
Let $\varrho:=\bar{F}_{g-h-1}^{\rm even} F_{h-1}^{\rm even}$ and 
$\varrho':=V_h\, (2\leq h\leq g-2)$ 
be the relators of $\mathcal{M}_g\, (g\geq 4)$ constructed in Section 4 of \cite{EN} and 
$\varsigma:=Q$ the positive relator of $\mathcal{M}_g\, (g\geq 2)$ 
constructed in Section 4 of \cite{endo}. 
Since $\varrho'$ is obtained by applying an $L$-substitution to $\varrho$, 
it turns out from Theorem \ref{bd} that $M_{\varrho'}$ is a rational blowdown of $M_{\varrho}$ 
along a copy of $C_2$. The Euler characteristic and the signature 
of $M_{\varrho'}$ are equal to $12g^2+6g+8gh-8h^2+7$ and $-6g^2-8g-4gh+4h^2-3$, respectively. 
$M_{\varsigma}-\nu F$ is simply-connected and $\varsigma$ includes one separating curve. 
The Euler characteristic and the signature of $M_{\varsigma}$ are 
$2g^2+7$ and $-(g^2+2g+3)$ for even $g$, and $2g^2+4g+7$ and $-(g+2)^2$ for odd $g$, 
respectively. 
We set $M_1:=M_{\varrho}\#_F M_{\varsigma}$ and $M_2:=M_{\varrho'}\#_F M_{\varsigma}$. 
It follows from Theorem 
\ref{exotic} that $M_1$ is homeomorphic but not diffeomorphic to $M_2\#\overline{\Bbb{CP}}^2$ 
and both of these do not dissolve. 
If we use $Q^n\, (n\geq 2)$ instead of $Q$, we obtain infinitely many pairs of 
homeomorphic but non-diffeomorphic $4$-manifolds for a fixed $g\, (\geq 4)$. 
\end{exmp}

\begin{exmp} Let $X_3$ and $X_{3,3}$ be the Lefschetz fibrations of genus $3$ defined in \S 4 of 
\cite{EG}. A positive relator $\varrho$ (resp. $\varrho'$) representing the monodromy of 
$X_3$ (resp. $X_{3,3}$) is given as follows (see Figure 3, Figure 4, and Figure 2 of \cite{EG}). 
\[ 
\varrho:=(c_1c_2x_1c_3rc_8c_8c_4x_2c_5c_6c_7)^3, \quad 
\varrho':=(\bar{y}_1x_1tvs_2c_8f_1c_8s_2\bar{x}_2r_3)^3, 
\]
where we put $r:={}_{f_1^{-1}}(c_4)$. 
We apply elementary transformations and simultaneous conjugations to 
$\varrho$ as follows. 
{\allowdisplaybreaks %
\begin{align*}
\varrho & =(c_1c_2x_1c_3rc_8c_8c_4x_2c_5c_6c_7)^3 
=(c_1c_2x_1c_3\cdot {}_{f_1^{-1}}(c_4)\cdot c_8c_8c_4x_2c_5c_6c_7)^3 \\
& \sim \; {}_{f_1^{-1}}(c_1c_2x_1c_3c_4c_8c_8\cdot {}_{f_1}(c_4)\cdot x_2c_5c_6c_7)^3 
\sim \; (c_1c_2x_1c_3c_4c_8c_8\cdot {}_{f_1}(c_4)\cdot x_2c_5c_6c_7)^3 \\
& \sim \; ({}_{c_1}(c_2)\cdot c_1x_1c_3c_4c_8c_8\cdot {}_{f_1}(c_4)\cdot x_2c_5c_6c_7)^3 \\
& \sim \; ({}_{c_1}(c_2)\cdot x_1c_3c_4c_8c_8\cdot {}_{f_1}(c_4)\cdot x_2c_5c_6c_7\cdot c_1)^3 \\
& \sim \; (c_1\cdot {}_{c_1}(c_2)\cdot x_1c_3c_4c_8c_8\cdot {}_{f_1}(c_4)\cdot x_2c_5c_6c_7)^3 \\
& \sim \; (c_1\cdot {}_{c_1}(c_2)\cdot x_1c_3c_4c_8c_8
\cdot c_5\cdot {}_{c_5^{-1}f_1}(c_4)\cdot {}_{c_5^{-1}}(x_2)\cdot c_6c_7)^3 \\
& \sim \; (c_1\cdot {}_{c_1}(c_2)\cdot x_1c_3c_4c_5c_8c_8
\cdot {}_{f_1c_5^{-1}}(c_4)\cdot \bar{x}_2\cdot c_6c_7)^3 \quad (\bar{x}_2:={}_{c_5^{-1}}(x_2)) \\
& \sim \; (c_1\cdot {}_{c_1}(c_2)\cdot x_1c_3c_5\cdot {}_{c_5^{-1}}(c_4)\cdot c_8c_8
\cdot {}_{f_1c_5^{-1}}(c_4)\cdot \bar{x}_2c_7\cdot {}_{c_7^{-1}}(c_6))^3  \\
& \sim \; (c_1\cdot {}_{c_1}(c_2)\cdot x_1c_3c_5\cdot s_2\cdot c_8c_8
\cdot {}_{f_1}(s_2)\cdot \bar{x}_2c_7\cdot r_3)^3  \quad 
(s_2:={}_{c_5^{-1}}(c_4),\, r_3:={}_{c_7^{-1}}(c_6)) \\
& \sim \; (c_1\cdot {}_{c_1}(c_2)\cdot x_1c_3c_5c_7s_2c_8c_8
\cdot {}_{f_1}(s_2)\cdot \bar{x}_2r_3)^3  \\
& \sim \; ({}_{c_1^2}(c_2)\cdot c_1\cdot x_1c_3c_5c_7s_2c_8c_8
\cdot {}_{f_1}(s_2)\cdot \bar{x}_2r_3)^3  \\
& \sim \; (\bar{y}_1x_1\cdot c_1c_3c_5c_7\cdot s_2c_8c_8
\cdot {}_{f_1}(s_2)\cdot \bar{x}_2r_3)^3=:\tau  \quad (\bar{y}_1:={}_{c_1^2}(c_2)) 
\end{align*}}
We apply elementary transformations to $\varrho'$ as follows. 
{\allowdisplaybreaks %
\begin{align*}
\varrho' & =\; (\bar{y}_1x_1tvs_2c_8f_1c_8s_2\bar{x}_2r_3)^3
\sim \; (\bar{y}_1x_1tvs_2c_8\cdot {}_{f_1}(c_8s_2\bar{x}_2r_3)\cdot f_1)^3 \\
& \sim \; (f_1\cdot \bar{y}_1x_1tvs_2c_8\cdot {}_{f_1}(c_8s_2\bar{x}_2r_3))^3 
\sim \; ({}_{f_1}(\bar{y}_1x_1)\cdot f_1tv\cdot s_2c_8\cdot {}_{f_1}(c_8s_2\bar{x}_2r_3))^3 \\
& = \; (\bar{y}_1x_1\cdot f_1tv\cdot s_2c_8c_8\cdot {}_{f_1}(s_2)\cdot\bar{x}_2r_3)^3 =:\tau'
\end{align*}}

\resizebox{!}{8.2cm}{\includegraphics{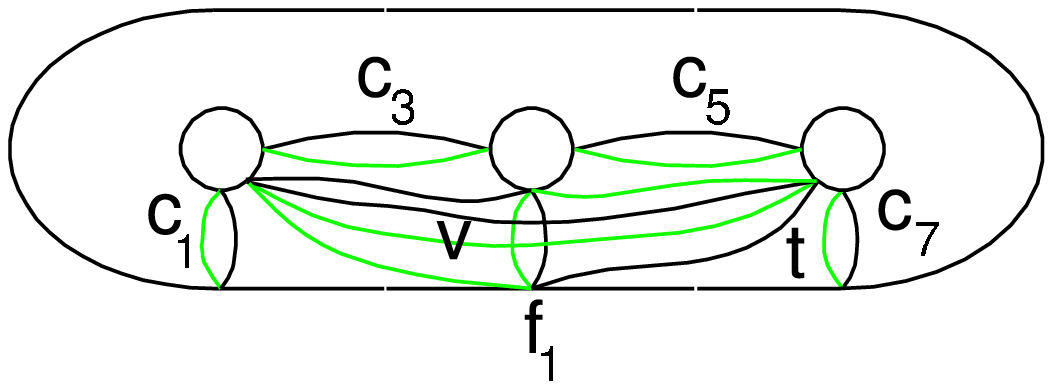}}

\vspace{-6cm}

\begin{center}
Figure 3
\end{center}

\resizebox{!}{9cm}{\includegraphics{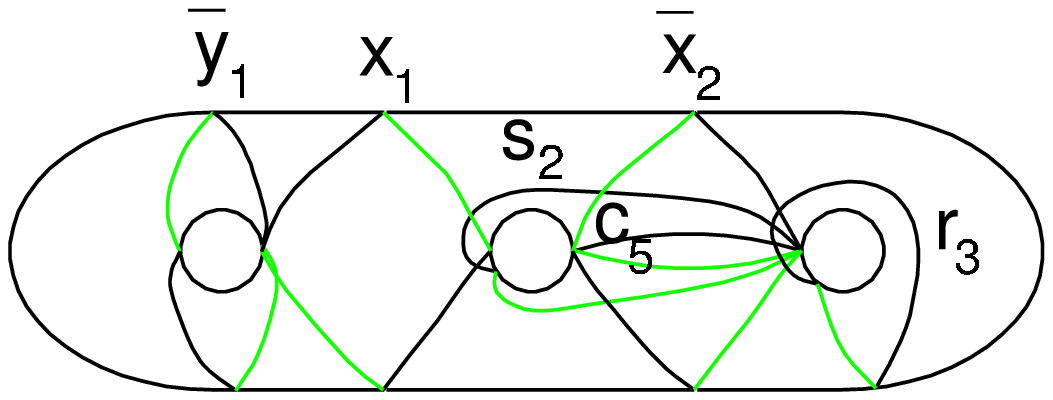}}

\vspace{-6.5cm}

\begin{center}
Figure 4
\end{center}

Thus $\tau'$ is obtained by applying three times $L$-substitutions to $\tau$ 
by virtue of the lantern relation $f_1tv=c_1c_3c_5c_7$, and 
$M_{\varrho'}=X_{3,3}$ turns out to be 
a rational blowdown of $M_{\varrho}=X_3$ along three copies of $C_2$ from Theorem \ref{bd}. 

We set $\varsigma:=(c_1c_2c_3c_4c_5c_6c_7^2c_6c_5c_4c_3c_2c_1)^2\in\mathcal{R}$ and 
put $M_1:=M_{\varrho}\#_F M_{\varsigma}$ and $M_2:=M_{\varrho'}\#_F M_{\varsigma}$. 
Both of $M_1$ and $M_2\# 3\overline{\Bbb{CP}}^2$ are simply-connected and  
have the Euler characteristic $56$ and signature $-36$. 
Hence Theorem \ref{exotic} tells us that 
$M_1$, $M_2\# 3\overline{\Bbb{CP}}^2$, and 
$\# 9\Bbb{CP}^2\# 45\overline{\Bbb{CP}}^2$ are homeomorphic but mutually 
non-diffeomorphic. 
\end{exmp}

We next exhibit an example of lantern substitution for genus $2$ fibrations 
and pose a problem about it. 

\begin{exmp}\label{g2}
Let $\varrho$ (resp. $\varrho'$) be a positive relator of $\mathcal{M}_2$ given as follows 
(see Figure 5, Figure 6, Figure 7, and Figure 4 of \cite{EG}). 
\begin{align*}
\varrho & :=(c_5c_4c_3c_2c_1^2c_2c_3c_4c_5)^2, \\
\varrho' & :={}_{c_3}(\delta) {}_{c_3c_4^{-1}}(x) \cdot \bar{k}\bar{h}c_5 \cdot {}_{c_1^{-4}}(c_2) 
{}_{c_1^{-1}}(c_2) {}_{c_2^{-1}c_3}(c_4) \cdot k\cdot {}_{c_2^{-1}}(h) {}_{c_2^{-1}c_3^{-1}}(c_4) 
{}_{c_1^2}(c_2) \cdot \bar{k}\bar{h}c_5c_4 
\end{align*}

\resizebox{!}{8.5cm}{\includegraphics{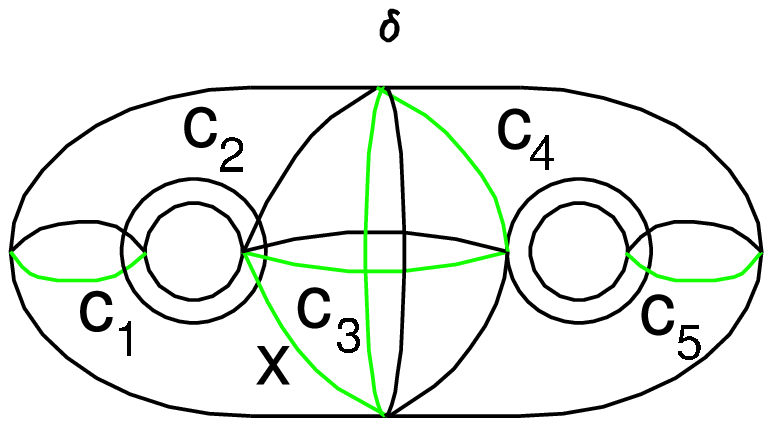}}

\vspace{-6cm}

\begin{center}
Figure 5
\end{center}

\vspace{5mm}

\resizebox{!}{7.3cm}{\includegraphics{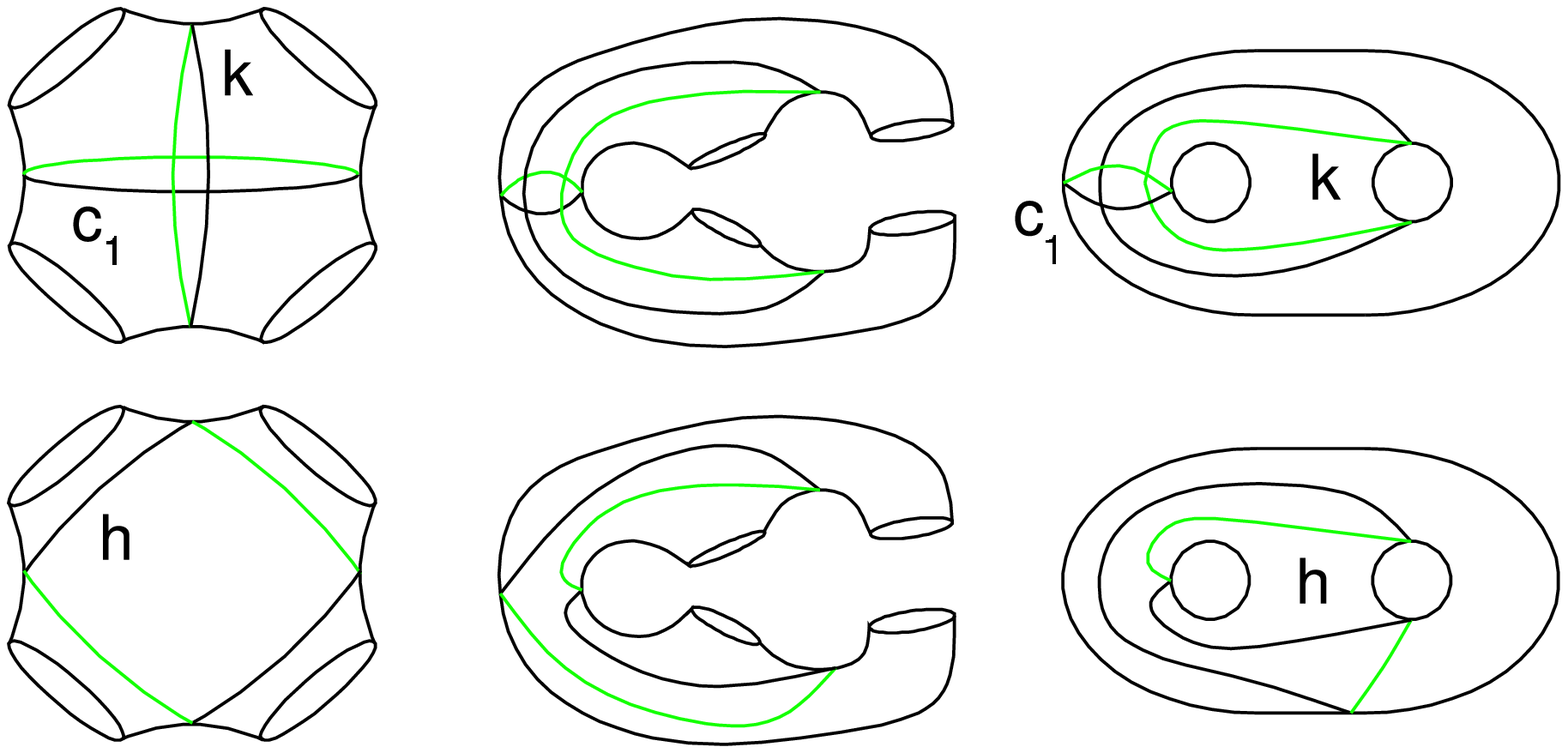}}

\vspace{-3cm}

\begin{center}
Figure 6
\end{center}

\resizebox{!}{10cm}{\includegraphics{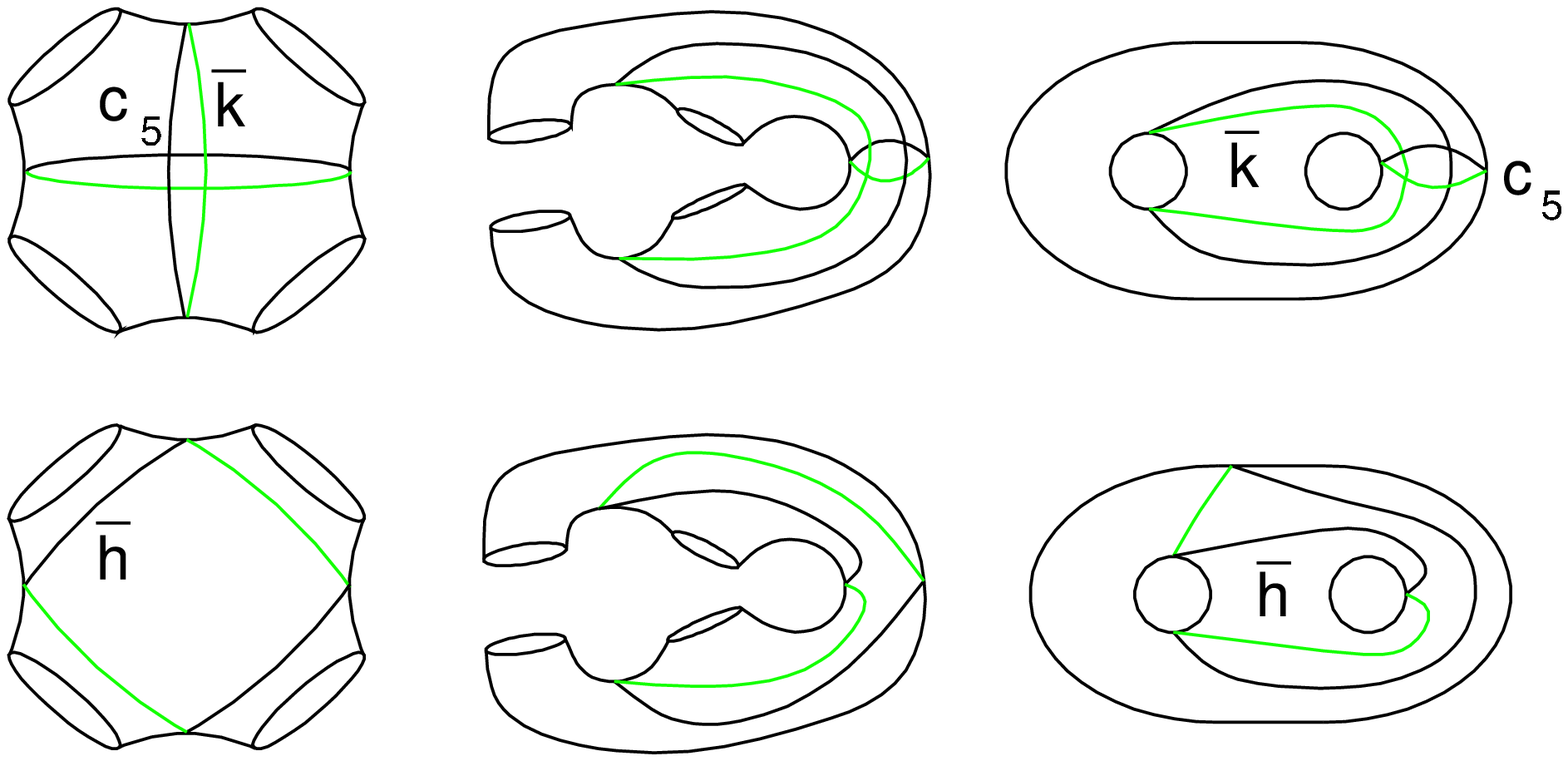}}

\vspace{-5.5cm}

\begin{center}
Figure 7
\end{center}

\medskip

Let $M_{\varrho}$ (resp. $M_{\varrho'}$) be the corresponding Lefschetz fibration 
of genus $2$ over $S^2$. 
It is well-known that $M_{\varrho}$ is diffeomorphic to 
$\Bbb{CP}^2\#13\overline{\Bbb{CP}}^2$ (cf. \cite{GS}). 
$\varrho'$ is obtained by applying elementary transformations 
and four times $L$-substitutions to $\varrho$ as follows. 
{\allowdisplaybreaks %
\begin{align*}
\varrho & =\; (c_5c_4c_3c_2c_1^2c_2c_3c_4c_5)^2 \\
& \sim \; c_5c_4c_3c_2c_1\cdot c_1c_2\cdot {}_{c_3}(c_4)\cdot c_3c_5\cdot
c_5c_3\cdot {}_{c_3^{-1}}(c_4)\cdot c_2c_1\cdot c_1c_2c_3c_4c_5 \\
& \overset{L}{\rightarrow} \; c_5c_4c_3c_2c_1\cdot c_1c_2\cdot {}_{c_3}(c_4)\cdot c_1kh 
\cdot {}_{c_3^{-1}}(c_4)\cdot c_2c_1\cdot c_1c_2c_3c_4c_5 \\
& \sim \; c_5c_4c_3\cdot c_1^2\cdot {}_{c_1^{-2}}(c_2)\cdot c_2\cdot {}_{c_3}(c_4)\cdot c_1kh 
\cdot {}_{c_3^{-1}}(c_4)\cdot c_2c_1\cdot c_1c_2c_3c_4c_5 \\
& \sim \; c_5^2c_1^2\cdot c_4c_3\cdot {}_{c_1^{-2}}(c_2)\cdot c_2\cdot {}_{c_3}(c_4)\cdot c_1kh 
\cdot {}_{c_3^{-1}}(c_4)\cdot c_2c_1\cdot c_1c_2c_3c_4 \\
& \overset{L}{\rightarrow} \; c_3\delta x\cdot c_4c_3\cdot {}_{c_1^{-2}}(c_2)\cdot 
c_2\cdot {}_{c_3}(c_4)\cdot c_1kh \cdot {}_{c_3^{-1}}(c_4)\cdot c_2c_1\cdot c_1c_2c_3c_4 \\
& \sim \; \delta x c_4c_3\cdot {}_{c_1^{-2}}(c_2)\cdot 
c_2\cdot {}_{c_3}(c_4)\cdot c_1kh \cdot {}_{c_3^{-1}}(c_4)\cdot c_2c_1\cdot c_1c_2c_3c_4c_3 \\
& \sim \; \delta x c_4c_3\cdot {}_{c_1^{-2}}(c_2)\cdot 
c_2\cdot {}_{c_3}(c_4)\cdot c_1kh \cdot {}_{c_3^{-1}}(c_4)\cdot c_2\cdot {}_{c_1^2}(c_2)
\cdot c_1^2c_3^2\cdot {}_{c_3^{-1}}(c_4) \\
& \overset{L}{\rightarrow} \; \delta x c_4c_3\cdot {}_{c_1^{-2}}(c_2)\cdot 
c_2\cdot {}_{c_3}(c_4)\cdot c_1kh \cdot {}_{c_3^{-1}}(c_4)\cdot c_2\cdot {}_{c_1^2}(c_2)
\cdot \bar{k}\bar{h}c_5\cdot {}_{c_3^{-1}}(c_4) \\
& \sim \; {}_{c_3^{-1}}(c_4)\cdot \delta x c_4c_3\cdot {}_{c_1^{-2}}(c_2)\cdot 
c_2c_1\cdot {}_{c_3}(c_4)\cdot kh \cdot {}_{c_3^{-1}}(c_4)\cdot c_2\cdot {}_{c_1^2}(c_2)
\cdot \bar{k}\bar{h}c_5 \\
& \sim \; {}_{c_3^{-1}}(c_4)\cdot \delta x c_4c_3\cdot {}_{c_1^{-2}}(c_2)\cdot 
c_2c_1c_2\cdot {}_{c_2^{-1}c_3}(c_4) \cdot k\cdot {}_{c_2^{-1}}(h) 
{}_{c_2^{-1}c_3^{-1}}(c_4) {}_{c_1^2}(c_2)\cdot \bar{k}\bar{h}c_5 \\
& \sim \; {}_{c_3^{-1}}(c_4)\cdot \delta x c_4c_3\cdot {}_{c_1^{-2}}(c_2)\cdot 
c_1c_2c_1\cdot {}_{c_2^{-1}c_3}(c_4) \cdot k\cdot {}_{c_2^{-1}}(h) 
{}_{c_2^{-1}c_3^{-1}}(c_4) {}_{c_1^2}(c_2)\cdot \bar{k}\bar{h}c_5 \\
& \sim \; {}_{c_3^{-1}}(c_4)\cdot \delta x c_4c_3c_1^2\cdot {}_{c_1^{-4}}(c_2) 
{}_{c_1^{-1}}(c_2) {}_{c_2^{-1}c_3}(c_4) \cdot k\cdot {}_{c_2^{-1}}(h) 
{}_{c_2^{-1}c_3^{-1}}(c_4) {}_{c_1^2}(c_2)\cdot \bar{k}\bar{h}c_5 \\
& \sim \; {}_{c_3^{-1}c_4c_3}(\delta) {}_{c_3^{-1}c_4c_3}(x) 
{}_{c_3^{-1}}(c_4)\cdot  c_4c_3c_1^2\cdot {}_{c_1^{-4}}(c_2) 
{}_{c_1^{-1}}(c_2) {}_{c_2^{-1}c_3}(c_4) \cdot k \\
& \qquad \qquad \cdot {}_{c_2^{-1}}(h) 
{}_{c_2^{-1}c_3^{-1}}(c_4) {}_{c_1^2}(c_2)\cdot \bar{k}\bar{h}c_5 \\
& \sim \; {}_{c_3^{-1}c_4c_3}(\delta) {}_{c_3^{-1}c_4c_3}(x) 
\cdot c_4c_3^2c_1^2\cdot {}_{c_1^{-4}}(c_2) 
{}_{c_1^{-1}}(c_2) {}_{c_2^{-1}c_3}(c_4) \cdot k \\
& \qquad \qquad \cdot {}_{c_2^{-1}}(h) 
{}_{c_2^{-1}c_3^{-1}}(c_4) {}_{c_1^2}(c_2)\cdot \bar{k}\bar{h}c_5 \\
& \sim \; c_4\cdot {}_{c_3}(\delta) {}_{c_3c_4^{-1}}(x) 
\cdot c_3^2c_1^2\cdot {}_{c_1^{-4}}(c_2) 
{}_{c_1^{-1}}(c_2) {}_{c_2^{-1}c_3}(c_4) \cdot k \cdot {}_{c_2^{-1}}(h) 
{}_{c_2^{-1}c_3^{-1}}(c_4) {}_{c_1^2}(c_2)\cdot \bar{k}\bar{h}c_5 \\
& \sim \; {}_{c_3}(\delta) {}_{c_3c_4^{-1}}(x) 
\cdot c_3^2c_1^2\cdot {}_{c_1^{-4}}(c_2) 
{}_{c_1^{-1}}(c_2) {}_{c_2^{-1}c_3}(c_4) \cdot k \cdot {}_{c_2^{-1}}(h) 
{}_{c_2^{-1}c_3^{-1}}(c_4) {}_{c_1^2}(c_2)\cdot \bar{k}\bar{h}c_5c_4 \\
& \overset{L}{\rightarrow} \; {}_{c_3}(\delta) {}_{c_3c_4^{-1}}(x) 
\cdot \bar{k}\bar{h}c_5\cdot {}_{c_1^{-4}}(c_2) 
{}_{c_1^{-1}}(c_2) {}_{c_2^{-1}c_3}(c_4) \cdot k \cdot {}_{c_2^{-1}}(h) 
{}_{c_2^{-1}c_3^{-1}}(c_4) {}_{c_1^2}(c_2)\cdot \bar{k}\bar{h}c_5c_4 \\
& = \; \varrho', 
\end{align*}}
where the symbol $\overset{L}{\rightarrow}$ stands for an $L$-substitution. 
Thus $M_{\varrho'}$ turns out to be a four times rational blowdown of 
$M_{\varrho}\approx \Bbb{CP}^2\#13\overline{\Bbb{CP}}^2$ 
along copies of $C_2$ from Theorem \ref{bd}. 

We set $\varsigma:=\varrho\in\mathcal{R}$ and 
put $M_1:=M_{\varrho}\#_F M_{\varsigma}$ and $M_2:=M_{\varrho'}\#_F M_{\varsigma}$. 
Both of $M_1$ and $M_2\# 4\overline{\Bbb{CP}}^2$ are simply-connected and  
have the Euler characteristic $36$ and signature $-24$. 
Hence Theorem \ref{exotic} tells us that $M_1$, $M_2\# 4\overline{\Bbb{CP}}^2$, and 
$\# 5\Bbb{CP}^2\# 29\overline{\Bbb{CP}}^2$ are homeomorphic but mutually 
non-diffeomorphic. 
\end{exmp}

We denote the manifold $M_{\varrho'}$ of Example \ref{g2} by $E$. 
Since $E$ is simply-connected and  
has the Euler characteristic $12$ and signature $-8$, $E$ is homeomorphic to 
$E(1)=\Bbb{CP}^2\# 9\overline{\Bbb{CP}}^2$ from Freedman's classification theorem. 

\begin{prob}
Does $E$ decompose into a non-trivial fiber sum of other Lefschetz fibrations? 
Is $E$ isomorphic to a fiber sum of two copies of Matsumoto's fibration? 
\end{prob}

If $E$ decomposes into a non-trivial fiber sum, 
then it is not diffeomorphic to $E(1)$ by virtue of Usher's theorem \cite{usher}. 
Matsumoto's fibration (Example B of \cite{yukiomat}) is a Lefschetz fibration of genus $2$ 
with $6$ non-separating, $2$ separating singular fibers, and its total space is 
diffeomorphic to $S^2\times T^2\# 4\overline{\Bbb{CP}}^2$. 
It is easy to see that an appropriately twisted fiber sum of two copies of 
Matsumoto's fibration is homeomorphic but not diffeomorphic to $E(1)$. 
Another possible way to examine the manifold $E$ would be to compute the Seiberg-Witten 
invariants of $E$ by the formula \cite{FS} of Fintushel and Stern. 

\section{Other relations}

We finally observe effects of substitutions for other relations. 
Luo \cite{luo} improved Gervais' infinite presentation \cite{gervais} of $\mathcal{M}_g$ to show 
that the set $\mathcal{R}$ of relators is normally generated by all commutativity, 
all braid, all $2$-chain, and all lantern relators. 
We briefly review definitions of these relations but lantern relation. 

Let $a,b$ be disjoint essential simple closed curves on $\Sigma_g$. 
The relation 
\[
t_at_b=t_bt_a
\]
in $\mathcal{M}_g$ is called a {\it commutativity relation}. 
A regular neighborhood $\Sigma$ of $a\cup b$ is the disjoint union of two annuli. 

Let $a,b$ be simple closed curves on $\Sigma_g$ which intersect transversely at one point. 
The relation 
\[
t_at_bt_a=t_bt_at_b
\]
in $\mathcal{M}_g$ is called a {\it braid relation}. 
A regular neighborhood $\Sigma$ of $a\cup b$ is a torus with one boundary component. 
Let $c$ be a simple closed curve parallel to the boundary of $\Sigma$. 
The relation 
\[
(t_at_b)^6=t_c
\]
in $\mathcal{M}_g$ is called a {\it chain relation of length 2}, or 
{\it 2-chain relation} in short. 

Both sides of each relation above correspond to Lefschetz fibrations over $D^2$ 
with fiber $\Sigma$. It is not difficult to draw Kirby diagrams of those manifolds and 
find out what they are (cf. \cite{GS}, Chapter 8). We actually obtain the following table. 

\bigskip 

\begin{tabular}{|c||c|c|c|} \hline 
relation & manifold for LHS & manifold for RHS & common boundary \\ \hline 
commutativity & $D^4\amalg D^4$ & $D^4\amalg D^4$ & $S^3\amalg S^3$ \\ \hline 
braid & $X(S^2,-2)$ & $X(S^2,-2)$ & $\Bbb{RP}^3$ \\ \hline 
$2$-chain & $M_c(2,3,6)$ & $X(T^2,-1)$ & $\Sigma(2,3,6)$ \\ \hline 
lantern & $C_2$ & $B_2$ & $L(4,1)$ \\ \hline 
\end{tabular} 

\bigskip 

\noindent
The symbol $X(B, e)$ stands for the total space of a $D^2$-bundle over $B$ with 
Euler number $e$. The Milnor fiber $M_c(2,3,6)$ and the Brieskorn manifold 
$\Sigma(2,3,6)$ are defined by 
\begin{align*}
M_c(2,3,6) & :=\{ (x,y,z)\in \Bbb{C}^3\, |\, x^2+y^3+z^6=\varepsilon \}\cap D^6, \\
\Sigma(2,3,6) & :=\{ (x,y,z)\in \Bbb{C}^3\, |\, x^2+y^3+z^6=0 \}\cap S^5
\end{align*}
(see \cite{GS}, Figure 8.13 for Kirby diagram). 
Substitutions for commutativity and braid relations do not change the original manifold 
(cf. 
\cite{yasui}, Figure 34 and 
\cite{endo}, Appendix A), 
whereas those for $2$-chain and lantern relations do. 

It might be interesting to extend the table above to that for various other relations such as 
chain relations of length $n\, (\geq 3)$, the star relation, and Matsumoto's relations 
\cite{matsumoto}. No relation seems to be known to correspond to a rational blowing down 
process along $C_p$ for $p\geq 3$.

\end{document}